\RequirePackage{fix-cm}
\documentclass[smallextended]{svjour3}       % onecolumn (second format)
\smartqed  % flush right qed marks, e.g. at end of proof
\usepackage{graphicx}
\begin{document}

\title{A novel and scalable Multigrid algorithm for many-core architectures}
%\subtitle{}


\author{Juli\'an Becerra-Sagredo \and Carlos M\'alaga \and Francisco Mandujano}


\institute{J. Becerra-Sagredo \at
              Inst. Ingenier\'ia Ed. 5 Cub. 307, UNAM, Ciudad Universitaria, 04510 M\'exico D.F., M\'exico \\
              Tel.: +525-562-24855\\
              \email{juliansagredo@gmail.com}           %  \\
           \and
           C. M\'alaga \and F. Mandujano \at
              Depto. F\'isica, Fac. Ciencias, UNAM, Ciudad Universitaria, 04510 M\'exico D.F., M\'exico
}

\date{Received: date / Accepted: date}

\maketitle

\begin{abstract}
Multigrid algorithms are among the fastest iterative methods known today for solving large 
linear and some non-linear systems of equations. Greatly optimized for serial operation, they 
still have a great potential for parallelism not fully realized.
In this work, we present a novel multigrid algorithm designed to work entirely inside many-core 
architectures like the graphics processing units (GPUs), without memory transfers between the 
GPU and the central processing unit (CPU), avoiding low bandwitdth communications. 
The algorithm is denoted as the high occupancy multigrid (HOMG) because it makes use of
entire grid operations with interpolations and relaxations fused into one task,
providing useful work for every thread in the grid. 
For a given accuracy, its number of operations scale linearly with the total number of nodes in 
the grid. Perfect scalability is observed for a large number of processors.
\end{abstract}

\section{Introduction}
\label{intro}
The multigrid algorithm \cite{Briggs,Trot} is one of the fastest methods for solving linear 
and non-linear systems of equations derived from a variety of problems, like 
numerical discretizations of partial differential equations and non-linear variational problems \cite{Hack,McCormick}. 
The main idea is to accelerate the convergence of an iterative method using a hierarchy of nested discretizations to  
perform resolution-dependent corrections that are passed between levels using interpolations. 
Multigrid has the main advantage over other methods that when used to solve the problems with a given accuray,
its number of operations often scale linearly with the number of discrete nodes used.
It has been widely studied and optimized for sequential operation during the last decades, 
providing very fast convergent algorithms 
that mostly lack fine grain parallelism because of the sequential use of updates during the same iteration.

In the last few years, hardware aspects have led to a paradigm change in the design of numerical methods. 
Performance 
improvements now come by parallelization and specialization and not any more by frequency scaling. 
Many-core parallel architectures, like GPUs, seem ideally suited for the 
acceleration of iterative methods like multigrid, given that its basic operations are essentially local and explicit.

Early works on multigrid for GPUs date back to 2003 \cite{Bolz,Goodnight}, when GPUs started 
to outperform CPUs  and control over the operations and memories of the GPU 
were available using programmable vertex shaders.
These implementations are many-core maps of a classic multigrid algorithm 
using recursive V-cycles. Additional interesting works in this direction are found in \cite{Strzodka,Thoman1,Grady,Sangild},
and more recently using Nvidia's compute unified device architecture (CUDA) 
 in \cite{Kazhdan,Feng,Elsen,Goddeke1,Rehman,Haase}. 
Nevertheless, a logical consequence of the use of a classic multigrid in 
many-core architectures is the appearance of performance penalties 
for coarse grids where the number of independent operations is reduced and memory latencies cannot be hided 
using multithreading. 

In order to avoid these penalties, 
a hybrid approach is preferred by some authors 
\cite{Thoman2,Goddeke2,Goddeke3} where the CPU 
is used for serial matrix inversions over coarse grids and the GPU for computing fine-grid relaxations.
This is a natural way to avoid the low parallelism 
of the coarse grids but with the penalty of the memory exchange 
between the random access memory (RAM) of the CPU and that of the GPU, 
which can be several orders of magnitude lower in bandwidth than the GPU fast off-chip memory.  

In this work, we modify the classic multigrid algorithm 
giving priority to keeping the data on the GPU, exploiting its fast memory levels and maintaining a constant 
occupancy. The result is a simple, fully parallel and perfectly scalable multigrid solver  
that works entirely inside the GPU without the need to communicate data to the CPU. 
We prove the perfect scalability of the algorithm using several GPUs with different numbers of processors.
The extension to multiple GPUs and clusters is beyond the scope of this work. 

\section{A high occupancy many-core multigrid (HOMG)}

Linear multigrid algorithms are iterative methods to solve large linear systems of the form
\begin{equation}
A {\bf u} = {\bf f},
\end{equation}
using a hierarchy of discretizations in nested grids  where level-dependent restrictions $P$, relaxations 
$R$,
and interpolations $S$ are performed in a variety of cycles. 

Classic V-cycles $V(\eta_1,\eta_2)$ consist of the recursive application of the following two-level correction scheme: 
(a)  $\eta_1$ pre-smoothing 
relaxations $v_h \leftarrow R(\eta_1) v_h$; (b) the restriction of the residual ${r_h} = { f_h} - A_h { v_h}$ 
to the next coarser grid $r_{2h} \leftarrow P r_h$; 
(c) the solution of the residual equation $A_{2h} {e_{2h}} = { r_{2h}}$; (d) the interpolation  
of the error to the fine grid and the correction of the solution $v_h \leftarrow v_h + S e_{2h}$; and (e) 
$\eta_2$ post-smoothing relaxations $v_h \leftarrow R(\eta_2) v_h$. 
Usually, interpolations are linear, restrictions are linearly weighted averages of neighboring sources  
and relaxations are one of the variations of Jacobi or Gauss-Sidel methods \cite{Briggs}.

A fundamental issue in many-core architectures like GPUs is the efficient use of 
its fast parallel fetches, optimized in every warp for contiguous memory locations and a large number of threads. 
Massive multithreading and the fast memory cache levels must be efficiently 
used to reduce cache misses and memory latencies. 
With this in mind, it is logical that some problems appear in multigrid algorithms 
as the grid is coarsened because the number of threads, 
proportional to the number of points in the grid, will be quite low, sometimes 
even lower than the number of processors available, and the processors will have very small amounts of work or will idle. 
Even colored Gauss-Sidel methods \cite{Briggs} could reduce the number of parallel threads and complicate the memory access
patterns. 
The coarsening of the grid and the consequent reduction of workload, basic properties of multigrid, 
become less attractive under the optic of GPUs, where operations over the entire grid are optimal and inexpensive
in one, two and three-dimensional structured grids. 

A multigrid algorithm that uses the entire grid for every iteration looks contradictory at first sight because 
it has more operations than the classic version, but 
it is attractive in many-core architectures because it keeps a constant optimal memory model fully exploiting 
its parallel fetching and processing capabilities.
One way to provide useful work for every thread 
during a coarse grid correction, is to fuse relaxations and 
interpolations in one task routine. 
Furthermore, we don't need to allocate different meshes because we could always act simultaneously 
in complementary subsets of the entire grid. 
One subset of points computes relaxations while
the rest perform interpolations of the previous iteration variations. 
This strategy provides a constant high occupancy of the cores. 

The damped Jacobi relaxation \cite{Briggs} is a natural choice for a fully parallel algorithm that only uses data produced 
in a previous iteration.
This method is equivalent to the explicit integration of a discretized version of 
\begin{equation}
\label{heat}
\frac{\partial u}{\partial t} = \kappa \left( \mathcal{L} (u) - f \right),
\end{equation}
where $\mathcal{L}(u)$ is any elliptic partial differential operator. 
For example, if $\mathcal{L} (u) = \nabla^{2} u$ and the Equation (\ref{heat}) is discretized using 
centered spatial differences, forward Euler integration in time
and nested meshes doubling the spatial discretization, the relaxation at level $m$ is given by
\begin{equation}
u_{i,j}^{n+1}  = u_{i,j}^{n} + \omega \left( (4u^{n}_{i,j} - u^{n}_{i-l,j} - u^{n}_{i+l,j} -u^{n}_{i,j-l} -u^{n}_{i,j+l}) - h_l^2  P_m(f_{i,j}) \right),
\end{equation}
where $l = 2^{m}$, $h_l = 2^{m} h$ and $P_m$ is the weighted average restriction to level $m$.

The high occupancy multigrid (HOMG) algorithm we propose is quite simple: 
(a) $\eta_1$ relaxations-interpolations on the coarsest level $v_{2^Mh} \leftarrow 
R(\eta_1) v_{2^Mh}$ and $v_h \leftarrow S v_{2^Mh}$; (b) $\eta_{2}$ 
relaxations-interpolations on the next finer level $v_{2^{(M-1)}h} \leftarrow 
R(\eta_2) v_{2^{(M-1)}h}$ and $v_h \leftarrow S v_{2^{(M-1)}h}$; (c) go back to a coarser level 
or descend to a finer level in a chosen pattern until reaching the finest level. 

\begin{figure}
\centering
\includegraphics*[width=7cm,height = 4cm]{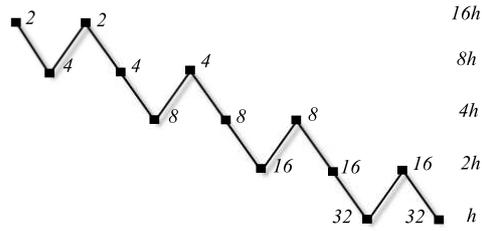}
\caption{Modified full multigrid cycle (MFMG). The number of iterations per level are doubling for every 
finer level until reaching a maximum number  ($MaxI$).}
\label{fig1}
\end{figure}

We have explored many sequences of coarse to fine grid levels and found the most efficient to be the 
modified full multigrid cycle (MFMG) presented in Fig. \ref{fig1}, with only two iterations in the coarsest level, increasingly 
doubling them for every finer level until reaching a chosen bound $MaxI$.

The restrictions are independent of the cycle and can be done in a pre-processing step using 
several independent arrays in memory. If memory savings are desirable, the restrictions can be done 
inside the cycle with just one extra memory array. 
The restrictions are performed inside the GPU using 
linear weighted averages organized as a  series of nested partial reductions over the entire mesh. 

After an MFMG cycle is complete, the residual equation is used to restart the cycles 
as many times as needed to reach any desired accuracy. 

\section{Power of two size grids}

\begin{figure}
\centering
\includegraphics*[width=4cm,height = 4cm]{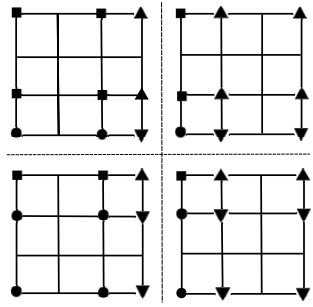}
\caption{Geometric strategy for multigrid in power of two grids. The grid is divided into four 
subdomains with a power of two number of points per side. Each set of symbols represents the location of the
coarse $2h$ grid in each subdomain. Some of them lie outside the subdomain. 
During a $2h$ relaxation-interpolation iteration, relaxations
are performed over all the points with a symbol and interpolations on the rest. 
The interpolations avoid cross-domain data.}
\label{fig2}
\end{figure}

Inherently to their architecture, GPUs perform optimally in grids with a power of two number of points per side
while multigrid algorithms are usually limited to nested grids.
We have designed a geometric strategy to allow the use of the HOMG in a two-dimensional grid 
with a power of two points per side. The procedure is illustrated in Fig. \ref{fig2}.

The grid is divided into four subdomains with a power of two nodes per side. The coarse grid points are marked 
for each subdomain relative to the corresponding external corner. 
These marked points would form a unique set in a nested configuration but here they form four different sets. 
The coarse subsets need to include central points, if not, some coarse modes, 
like the zero mode, can not be resolved. Due to the absence of central points we 
mark the coarse grid points that lie just outside the corresponding subdomain. 
We perform relaxations in all the marked points and interpolations on the 
rest of the fine grid. The interpolations make use only of points in the same subdomain to avoid cross-domain 
interpolations.

\section{Results}

We test our algorithm solving the two-dimensional Poisson equation 
\begin{equation}
\label{poisson}
\nabla^{2} u = f
\end{equation} 
on the unit square $\Omega = [0,1] \times [0,1]$ with Dirichlet boundary condition $u = 0$ on the 
edge $\partial \Omega$. 

We choose the polynomial analytic solution
 \begin{equation}
 u (x_1,x_2) = - x_1^2  x_2^2 (1-x_1^2) (1- x_2^2)
 \end{equation}
 corresponding to the source 
 \begin{equation}
 \label{analytic1}
 f(x_1,x_2) = -2(x_2^2(1- 6x_1^2)(1-x_2^2) +x_1^2( 1- 6x_2^2)(1-x_1^2)).
 \end{equation}

This solution has a wide spectrum with components in all the modes of the grid. 
This characteristic is quite useful to study the convergence of multigrid
algorithms like the HOMG, and has been previously used in \cite{Briggs}.

\begin{figure}
\centering
\includegraphics*[width=11.8cm,height = 4.8cm]{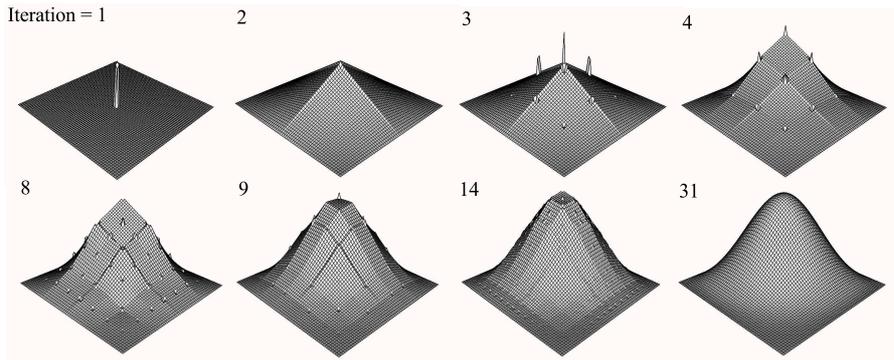}
\caption{High occupancy multigrid (HOMG) approximation to the solution of Eq. (\ref{poisson}) with the polynomial source (\ref{analytic1})
in a $128^2$ grid. 
The approximation is shown for different iterations of the modified full multigrid cycle (MFMG) 
with a maximum of two iterations per level ($MaxI = 2$).
}
\label{fig3}
\end{figure}

Fig. \ref{fig3} shows the sequential approximation of the solution with the HOMG algorithm 
in one MFMG cycle with a maximum of two iterations per level ($MaxI = 2$) and a $128^2$ grid. 
During the first iteration we see a spike because 
only the middle points perform a relaxation and the interpolations
have null sources. Remember that the interpolations are delayed one iteration because they are 
performed in parallel together with the relaxations. 
The second iteration shows that the relaxations have stalled while 
the interpolations form the zero mode approximation with pyramid shape. 
The third iteration shows clearly the next finer level relaxations-interpolations and so on. 
The MFMG cycle has a total of $31$ iterations for a $128^2$ grid. 

The error of the approximations is computed using the normalized $L_1$-error
\begin{equation}
Error = \frac{1}{\int_{\Omega} | u | d \Omega} \int_{\Omega} | u - v_h | d \Omega.
\end{equation}

\begin{figure}
\center
\includegraphics*[width=10cm,height = 6cm]{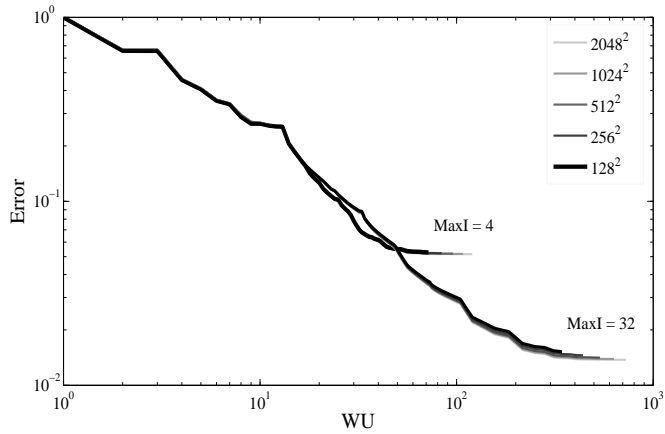} 
\caption{$L_1$-error (Error) convergence of the modified multigrid cycle for different domain 
sizes and maximum number of 
iterations per level ($MaxI$). The working units (WU) correspond to entire grid operations. }
\label{fig4}
\end{figure}

Fig. \ref{fig4} shows the convergence of the $L_1$-error for the HOMG method for one MFMG cycle and two 
different choices of the maximum iterations per cycle $MaxI$. 
The results show all the desirable properties of a multigrid algorithm. 
Cycles with the same $MaxI$ descend over the same convergence path for different resolutions. 
Higher resolutions continue the smoothing of the solution obtained with coarser grids. 
Increasing $MaxI$ allows the cycles to reach lower errors with the penalty of a large number of iterations 
per cycle. $MaxI = 4$ reaches an error of $5$\% and $MaxI = 32$ reaches an error of $1.6$\%. 
Operations over the entire grid are denoted as a working unit (WU).
Given a desired accuracy, the algorithm performs a number of operations linearly proportional 
to the points in the grid. Remember that after an MFMG cycle is complete, 
the residual equation is used to restart the cycles 
as many times as needed to reach any desired accuracy. 

\begin{figure}
\center
\includegraphics*[width=10cm,height = 6cm]{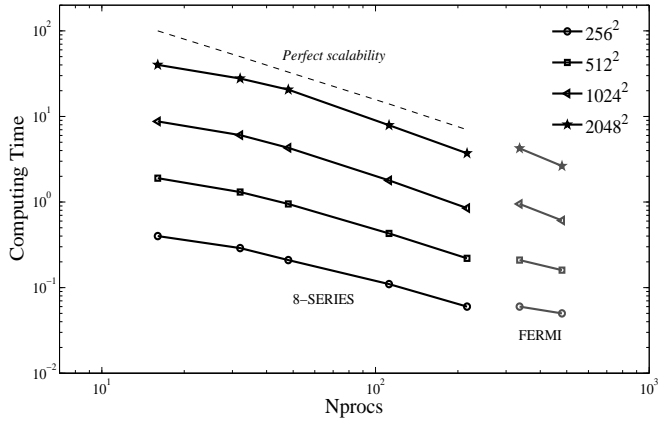} 
\caption{Computing time in seconds for  $MaxI = 32$ in GPUs with different architectures and number of 
processors (Nprocs). The perfect linear scalability curve is shown for comparisson. }
\label{fig5}
\end{figure}

We explore the scalability of the HOMG algorithm on GPUs with different architectures and number of processors.
The test consists of one MFMG cycle with $MaxI = 32$ and different resolutions.  
The GPUs used are all GeForce from Nvidia. We use C language for CUDA  as the 
application programming interface.
Five GPUs of the modified 8-series architecture: a 310M with 16 processors, a 9500GT with 32 processors, a 330M with 48 processors, a 9800GT with 112 processors and a GTX260 with 216 processors. 
And two GPUs of the GeForce Fermi architecture: a GTX460 with 336 processors and a GTX480 with 480 processors. 
All the computations are done in single precision. 

The scalability results are shown in Fig. \ref{fig5}. The slopes are moderate for a low number of processors and small grids. 
For a large number of processors and large grids the algorithm reaches perfect scalability.

\section{Conclusion}

In this article, we present a novel multigrid algorithm specially designed to work entirely inside many-core architectures 
like GPUs without memory transfers between the GPU and the CPU. The algorithm makes use of entire grid operations 
even for coarse grid corrections. Interpolations and relaxations are fused into one task giving useful work for every thread 
in the grid. In this way the algorithm has full multithreading and fix memory patterns, allowing the full exploitation of 
the fast memory models of the GPU, efficiently hiding cash misses and memory latencies. 

The algorithm is denoted as the high occupancy multigrid (HOMG) algorithm because multithreading 
and useful work per thread are kept constantly high. The algorithm is combined with a modified full multigrid cycle 
(MFMG) to reach a high efficiency. For a given accuracy, the operations of the HOMG scale linearly with the 
total number of nodes. Perfect scalability is observed for a large number of processors and large grids.

\end{document}